\documentclass{ifacconf}

\usepackage{graphicx}      % include this line if your document contains figures
\usepackage{natbib}        % required for bibliography

\usepackage{amssymb} %amsfonts, latexsym, mathtools}
\usepackage{amsmath}

\usepackage{enumerate}

\def\N{\mathbb{N}}
\def\R{\mathbb{R}}
\def\KK{\mathcal{K}}
\def\OO{\mathcal{O}}
\def\eps{\varepsilon}

\newtheorem{theorem}{Theorem}[section]
\newtheorem{definition}[theorem]{Definition}
\newtheorem{proposition}[theorem]{Proposition}

\begin{document}
\begin{frontmatter}

\title{Overcoming the curse of dimensionality\\ for approximating Lyapunov functions\\ with deep neural networks\\ under a small-gain condition} 

%\thanks[footnoteinfo]{place thanks here}

\author[First]{Lars Gr\"{u}ne}

\address[First]{Mathematical Institute, 
	University of Bayreuth, Bayreuth, Germany, e-mail: lars.gruene@uni-bayreuth.de}

\begin{abstract}                % Abstract of not more than 250 words.
We propose a deep neural network architecture for storing approximate Lyapunov functions of systems of ordinary differential equations. Under a small-gain condition on the system, the number of neurons needed for an approximation of a Lyapunov function with fixed accuracy grows only polynomially in the state dimension, i.e., the proposed approach is able to overcome the curse of dimensionality.
\end{abstract}


\end{frontmatter}
%===============================================================================

\section{Introduction}
The computation of Lyapunov functions has attracted significant attention during the last decades. The reason is that the knowledge of a Lyapunov function not only serves as a certificate for asymptotic stability of an equilibrium but also allows to give estimates about its domain of attraction or to quantify its robustness with respect to perturbations, for instance, in the sense of input-to-state stability. As explicit analytic expressions for Lyapunov functions are typically not available, numerical methods need to be used.

Known approaches use, e.g., representations by radial basis functions as in \cite{Gies07}, piecewise affine functions, see \cite{Hafs07}, or sum-of-squares techniques, see \cite{AndP15} and the references therein. For a comprehensive overview we refer to the survey by \cite{GieH15}.

The usual approaches have in common that the number of degrees of freedom needed for storing the Lyapunov function (or an approximation thereof with a fixed approximation error) grows very rapidly --- typically exponentially --- with the dimension of the state space. This is the well known curse of dimensionality, which leads to the fact that the mentioned approaches are confined to low dimensional systems.

In general, the same is true if a deep neural network is used as an approximation architecture. While it is known that such a network can approximate every $C^1$-function arbitrarily well, see \cite{Cybe89,HoSW89}, the number of neurons needed for this purpose typically grows exponentially with the state dimension, see \cite[Theorem 2.1]{Mhas96} or Theorem \ref{thm:univ}, below. However, this situation changes drastically
if additional structural assumptions are imposed, which is the approach we follow in this paper. 

More precisely, we show that under a small-gain condition a deep neural network can store an approximation of a Lyapunov function with a given required accuracy using a number of neurons that grows only polynomially with the dimension of the system. In other words, under a small-gain condition the curse of dimensionality can be avoided, at least in terms of storage requirement. 

The study in this paper was inspired by a couple of recent similar results for solutions of partial differential equations and optimal value functions, such as \cite{DaLM19,HuJT19,ReiZ19}, and our proof technique makes use of some arguments from \cite{PMRML17}. However, the --- to the best of our knowledge --- entirely new insight is that a small-gain condition leads to a class of Lyapunov functions that allows for applying these techniques. This is, in our opinion, the main contribution of this paper.

\section{Problem Formulation}

We consider nonlinear ordinary differential equations of the form
\begin{equation} \dot x(t) = f(x(t)) \label{eq:sys}\end{equation}

with a Lipschitz continuous $f:\R^n\to\R^n$. We assume that $x=0$ is a globally\footnote{In order to avoid technicalities in the exposition we consider global asymptotic stability here. The subsequent arguments could be adapted if the domain of attraction of $x=0$ is not the whole $\R^n$.}asymptotically stable equilibrium of~\eqref{eq:sys}. 

It is well known (see, e.g., \cite{Hahn67}) that this is the case if and only if there exists a $C^1$-function $V:\R^n\to\R$ satisfying $V(x)\ge 0$ for all $x\in\R^n$, $V(x)>0$ for all $x\ne 0$, $V(x)\to \infty$ as $\|x\|\to\infty$ and 
\begin{equation} DV(x)f(x) \le - h(x) \label{eq:lf}\end{equation}

for a function $h:\R^n\to\R$ with $h(x)>0$ for all $x\ne 0$. It is our goal in this paper to design a neural network that is able to store an approximation of such a Lyapunov function on a compact set $K\subset\R^n$ in an efficient manner. To this end, we exploit the particular structure of $V$ that is induced by a small-gain condition.

\section{Small-gain conditions and Lyapunov functions}

In small-gain theory, the system \eqref{eq:sys} is divided into $s$ subsystems $\Sigma_i$ of dimensions $d_i$, $i=1,\ldots,s$. To this end, the state vector $x=(x_1,\ldots,x_n)^T$ and the vector field $f$ are split up as
\[ x = \left(\begin{array}{c} z_1\\ z_2\\ \vdots \\ z_s\end{array}\right) \; \mbox{ and } \; f(x) = \left(\begin{array}{c} f_1(x)\\ f_2(x)\\ \vdots \\ f_s(x)\end{array}\right),\]
with $z_i\in\R^{d_i}$ and $f_i:\R^n\to\R^{d_i}$ denoting the state and dynamics of each $\Sigma_i$, $i=1,\ldots,s$. With 
\[ z_{-i} := \left(\begin{array}{c} z_1\\ \vdots \\ z_{i-1}\\z_{i+1}\\ \vdots \\ z_s\end{array}\right) \]
and by rearranging the arguments of the $f_i$, the dynamics of each $\Sigma_i$ can then be written as
\[ \dot z_i(t) = f_i(z_i(t),z_{-i}(t)), \quad i=1,\ldots,s. \]
Nonlinear small gain theory that we apply here relies on the input-to-state stability property introduced in \cite{Sont89}. It goes back to \cite{JiTP94} and in the form for large-scale systems we require here it was developed in the thesis by \cite{Ruef07} and in a series of papers around 2010, see, e.g., \cite{DaRW10,DaIW11} and the references therein. ISS small-gain conditions can be based on trajectories or Lyapunov functions and exist in different variants. Here, we use the variant that is most convenient for obtaining approximation results because it yields a smooth Lyapunov function. We briefly discuss one other variant in Section \ref{sec:discussion}\eqref{it:maxlf}.

For formulating the small gain condition, we assume that for the subsystems $\Sigma_i$ there exist $C^1$ ISS-Lyapunov functions $V_i:\R^{d_i}\to\R$, satisfying for all $z_i\in\R^{d_i}$ $z_{-i}\in\R^{n-d_i}$
\[ DV_i(z_i)f_i(z_i,z_{-i}) \le -\alpha_{i}(z_i) + \sum_{j\ne i} \gamma_{ij}(V_j(z_j)) \]
with functions $\alpha_i$, $\gamma_{ij}\in\KK_\infty$,\footnote{As usual, we define $\KK_\infty$ to be the space of continuous functions $\alpha:[0,\infty)\to[0,\infty)$ with $\alpha(0)=0$ and $\alpha$ is strictly increasing to $\infty$.} $i,j=1\ldots,s$, $i\ne j$. Setting $\gamma_{ii}:= 0$, we define the map $\Gamma:[0,\infty)^s\to[0,\infty)^s$ by 
\[ \Gamma(r) := \left(\sum_{j=1}^s\gamma_{1j}(r_j),\ldots,\sum_{j=1}^s\gamma_{sj}(r_j)\right)^T\]
and the diagonal operator $A:[0,\infty)^s\to[0,\infty)^s$ by 
\[ A(r) := \left(\alpha_{1}(r_1),\ldots,\alpha_{s}(r_s)\right)^T.\]

\begin{definition} \label{def:sg} We say that \eqref{eq:sys} satisfies the small-gain condition, if there is a decomposition into subsystems $\Sigma_i$, $i=1,\ldots,s$, with ISS Lyapunov functions $V_i$ satisfying the following condition: there are bounded positive definite\footnote{A continuous function $\rho:[0,\infty)\to[0,\infty)$ is called positive definite if $\rho(0)=0$ and $\rho(r)>0$ for all $r>0$.} functions $\eta_i$, $i=1,\ldots,s$, satisfying $\int_0^\infty \eta_i(\alpha_i(r))dr = \infty$ and such that for $\eta=(\eta_1,\ldots,\eta_s)^T$ the inequality
\[ \eta(r)^T\Gamma\circ A(r) < \eta(r)^Tr\]
holds for all $r\in[0,\infty)^s$ with $r\ne 0$.
\end{definition}

The following theorem then follows from Theorem 4.1 in \cite{DaIW11}.

\begin{theorem}
Assume that the small-gain conditions from Definition \ref{def:sg} hold. Then a Lyapunov function for \eqref{eq:sys} is given by 
\[ V(x) = \sum_{i=1}^s \widehat V_i(z_i), \]
for the $C^1$-functions $\widehat V_i:\R^{d_i}\to\R$ given by
\[ \widehat V_i(z_i):= \int_0^{V_i(z_i)} \lambda_i(\tau)d\tau \]
where $\lambda_i(\tau) = \eta_i(\alpha_i(\tau))$. 
\label{thm:sg}\end{theorem}

In \cite{DaIW11}, the property from Definition~\ref{def:sg} is called a \emph{weak} small-gain condition. This is because if the system \eqref{eq:sys} has an additional input (that is taken into account in the assumptions on the $V_i$), then the construction of $V$ yields an integral ISS Lyapunov function as opposed to an ISS Lyapunov function. Under a stronger version of the small-gain condition, the same construction yields an ISS Lyapunov function. We briefly discuss corresponding extensions of our approach in Section \ref{sec:discussion}\eqref{it:u}.

\section{Deep neural networks}

\begin{figure}[htb]
\begin{center}
\includegraphics[width=8cm]{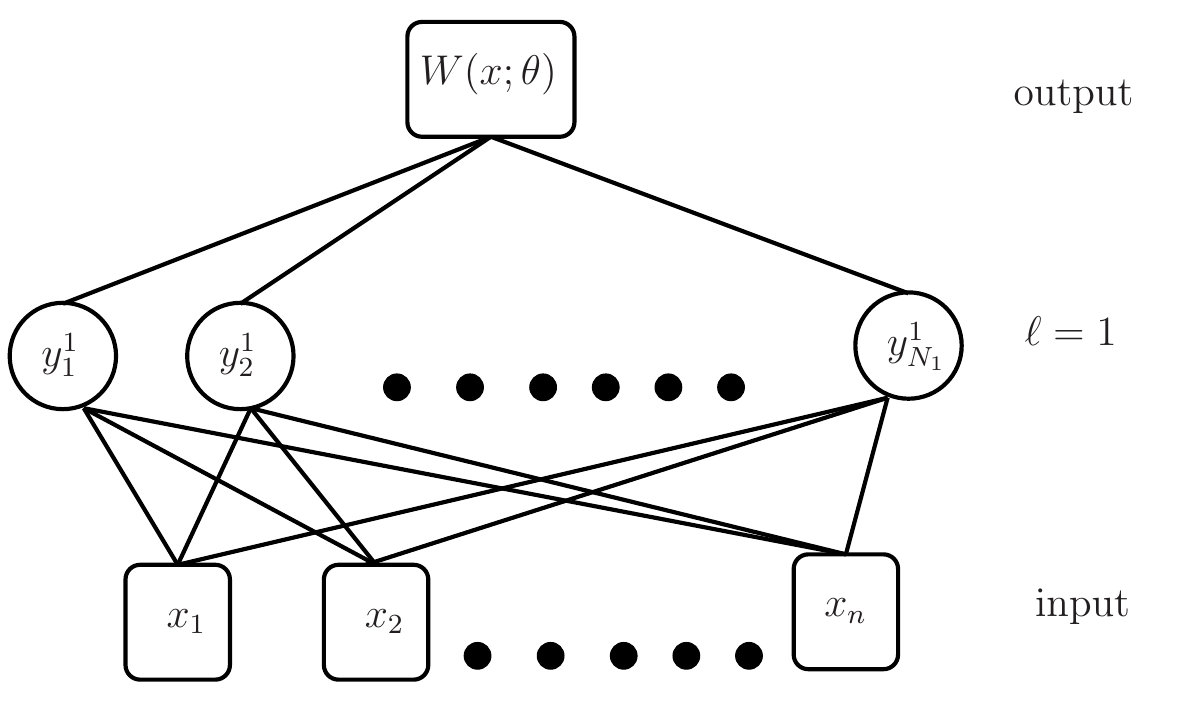}\\[10mm]
\includegraphics[width=8cm]{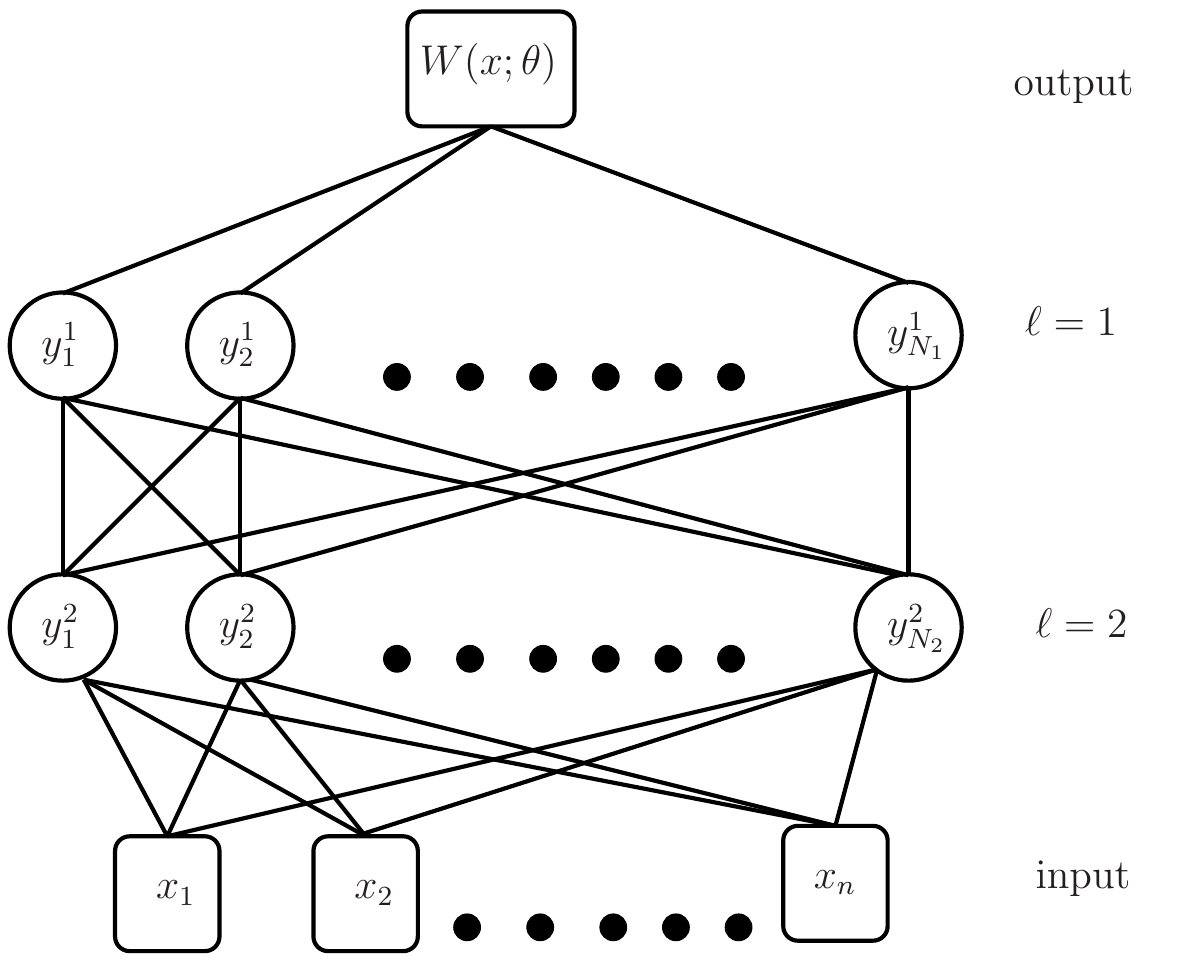}
\end{center}
\caption{Neural network with $1$ and $2$ hidden layers}
\label{fig:gen_nn}
\end{figure}

A deep neural network is a computational architecture that has several inputs, which are processed through $\ell\ge 1$ hidden layers of neurons. The values in the neurons of the highest layer ($\ell=1$) is then used in order to compute the output of the network. In this paper, we will only consider feedforward networks, in which the input is processed consecutively through the layers $\ell$, $\ell-1$, \ldots, $1$. For our purpose of representing Lyapunov functions, we will use networks with the input vector $x=(x_1,\ldots,x_n)^T\in\R^n$ and one output $W(x;\theta)\in\R$. Here, the vector $\theta\in\R^P$ represents the free parameters in the network that need to be tuned (or ``learned'') in order to obtain the desired output. In our case, the output shall approximate the Lyapunov function, i.e., we want to choose $\theta$ such that $W(x;\theta)\approx V(x)$ for all $x\in K$, where $K$ is a fixed compact set containing the origin. 
Figure \ref{fig:gen_nn} shows generic neural networks with one and two hidden layers. 

Here, the topmost layer is the output layer, followed by one or two hidden layers numbered with $\ell$, and the input layer. 
The number $\ell_{\max}$ determines the number of hidden layers, here $\ell_{\max}=1$ or $2$. Each hidden layer consists of $N_\ell$ neurons and the overall number of neurons in the hidden layers is denoted by $N=\sum_{\ell=1}^{\ell_{\max}} N_\ell$. The neurons are indexed using the number of their layer $\ell$ and their position in the layer $k$. Every neuron has a value $y^\ell_k$, for each layer these values are collected in the vector $y^{\ell}=(y^\ell_1,\ldots,y^\ell_{N_\ell})^T$. The values of the neurons at the lowest level are given by the inputs, i.e., $y^{\ell_{\max}+1}=x$. The values of the neurons in the hidden layers are determined by the formula
\[ y_k^\ell = \sigma^\ell(w_k^\ell \cdot y^{\ell+1} + b_k^\ell),\] 
for $k=1,\ldots N_\ell$, where $\sigma^\ell:\R\to\R$ is a so called activation function and $w_k^\ell\in\R^{N_{\ell+1}}$, $a_k^\ell,b_k^\ell\in\R$ are the parameters of the layer. With $x\cdot y$ we denote the Euclidean scalar product between two vectors $x,y\in\R^n$. In the output layer, the values from the topmost hidden layer $\ell=1$ are affine linearly combined to deliver the output, i.e., 
\[ W(x;\theta) = \sum_{k=1}^{N_1} y_k^1 + c = \sum_{k=1}^{N_1} a_k \sigma^1(w_k^1 \cdot y^2 + b_k^1) + c.\]
The vector $\theta$ collects all parameters $a_k$, $w_k^\ell$, $b_k^\ell$, $c$ of the network. 

In case of one hidden layer, in which $\ell_{\max}=1$ and thus $y^2=y^{\ell_{\max}+1}=x$, we obtain the closed-form expresssion
\[ W(x;\theta) = \sum_{k=1}^{N_1} a_k^1 \sigma^1(w_k^1\cdot x + b_k^1) + c.\]
The universal approximation theorem states that a neural network with one hidden layer can approximate all smooth functions arbitrarily well. In its qualitative version, going back to \cite{Cybe89,HoSW89}, it states that the set of functions that can be approximated by neural networks with one hidden layer is dense in the set of continuous functions. In Theorem \ref{thm:univ}, below, we state a quantitative version, given as Theorem 1 in \cite{PMRML17}, which is a reformulation of Theorem 2.1 in \cite{Mhas96}. 

For its formulation we fix $C>0$ and compact sets $K_n\subseteq [-C,C]^n\subset\R^n$. For a continuous function $g:K_n\to \R$ we define the infinity-norm over $K$ as
\[ \|g\|_{\infty,K_n} := \max_{x\in K_n}|g(x)|. \]
Given a fixed constant $C>0$, we then define the set of functions 
\[ W_m^n := \left\{ g\in C^m(K_n,\R) \,\left|\, \sum_{1\le |\alpha|\le m} \| D_\alpha g\|_{\infty,K} \le 1\right.\right\} \]
where $C^m(K_n,\R)$ denoted the functions from $K_n$ to $\R$ that are $m$-times continuously differentiable, $\alpha$ are multiindices of length $|\alpha|$ with entries $\alpha_i\in\{1,\ldots,n\}$, $i=1,\ldots,|\alpha|$ and $D_\alpha g = \partial g^{|\alpha|}/\partial \alpha_1\ldots \partial\alpha_{|\alpha|}$ denotes the $m$-th directional derivative with respect to $\alpha$.
%In \cite{PMRML17} the set $W_m^n$ is defined with $C=1$, however, the results are easily transferred to arbitrary $C>0$ by a straightforward rescaling. 

\begin{theorem} Let $\sigma:\R\to\R$ be infinitely differentiable and not a polynomial. Then, for any $\eps>0$, a neural network with one hidden layer provides an approximation $\inf_{\theta\in\R^P} \|W(x;\theta) - g(x)\|_{\infty,K_n} \le \eps$ for all $g\in W_m^n$ with a number of $N$ of neurons satisfying
\[ N = \OO\left(\eps^{-\frac{n}{m}}\right)\]
and this is the best possible.
\label{thm:univ}\end{theorem}

Theorem \ref{thm:univ} implies that one can readily use a network with one hidden layer for approximating Lyapunov functions. However, in general the number $N$ of neurons needed for a fixed approximation accuracy $\eps>0$ grows exponentially in $n$, and so does the number of parameters in $\theta$. This means that the storage requirement as well as the effort to determine $\theta$ easily exeeds all reasonable bounds already for moderate dimensions $n$. Hence, this approach also suffers from the curse of dimensionality. In the next section, we will therefore exploit the particular structure of small-gain based Lyapunov functions in order to obtain neural networks with (asymptotically) much lower $N$.

\section{Results}

\subsection{The case of known subsystems}\label{sec:known}

%We develop our results for systems \eqref{eq:sys} with arbitrary dimension $n$ under small type structural assumptions. 
For our first result, for fixed $d^{\max}\in\N$ we consider the family $F_1^{d_{\max}}$ of Lipschitz maps $f:\R^n\to\R^n$, $n\in\N$, for which \eqref{eq:sys} satisfies the small-gain condition from Definition~\ref{def:sg} with $\max_{i=1,\ldots,s} d_i \le d_{\max}$. We assume that for each $f$ in this family we know the subsystems $\Sigma_i$ satisfying the small-gain condition from Definition \ref{def:sg}, i.e., their dimensions $d_i$, $i=1,\ldots,s$ and states $z_i$.

\begin{figure}[htb]
\begin{center}
\includegraphics[width=8cm]{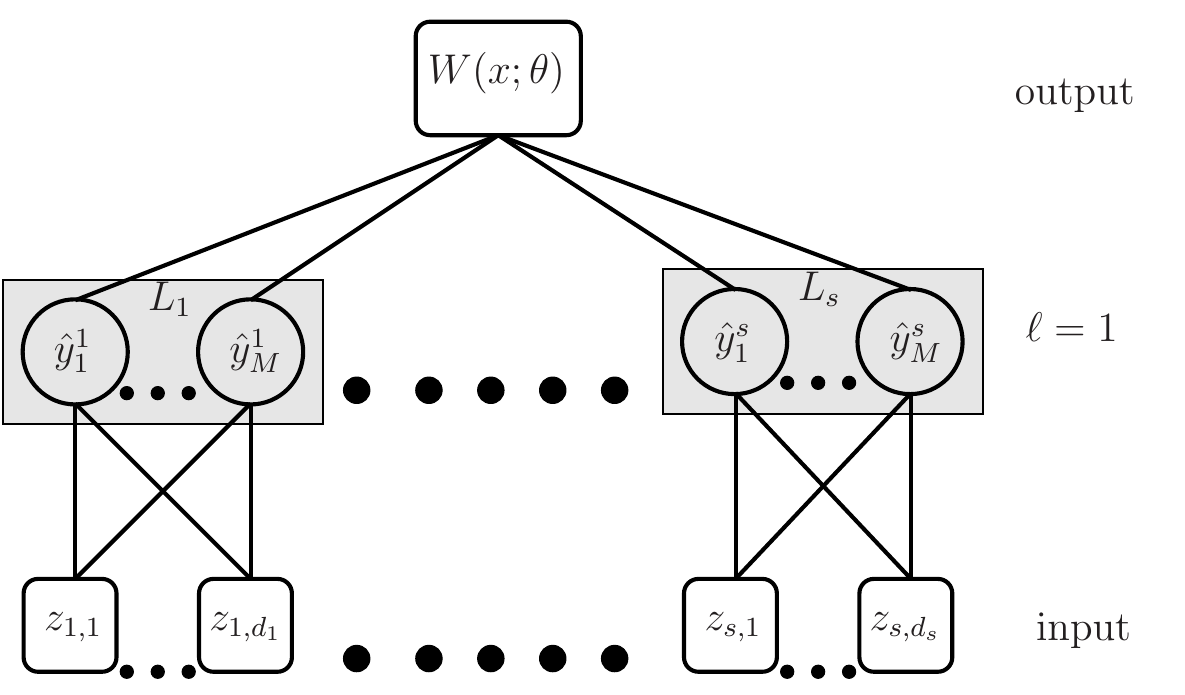}
\end{center}
\caption{Neural network for Lyapunov functions, $f\in F_1^{d_{\max}}$}
\label{fig:Lf_nn_l1}
\end{figure}

For this situation, we use a network with one hidden layer of the form depicted in Figure \ref{fig:Lf_nn_l1}.

In this network, the (only) layer for $\ell=1$ consists of $s$ sublayers $L_1, \ldots, L_s$. The input of each of the neurons in $L_i$ is the state vector $z_i=(z_{i,1},\ldots,z_{i,d_i})^T$ of the subsystem $\Sigma_i$, which forms a part of the state vector $x$. We assume that every sublayer $L_i$ has $M$ neurons, whose parameters and values are denoted by, respectively, $\hat w_k^i$, $\hat a_k^i$, $\hat b_k^i$ and $\hat y_k^i$, $k=1,\ldots,d_i$. Since $s\le n$, the layer contains $N^1=sM\le nM$ neurons, which is also the total number $N$ of neurons in the hidden layers. The values $\hat y_k^i$ are then given by 
\[ \hat y_k^i = \sigma^1(\hat w_k^i \cdot z_i + \hat b_k^i)\] 
and the overall output of the network is 
\[ W(x;\theta) = \sum_{i=1}^s \sum_{k=1}^{d_i} \hat a_k^i\sigma^1(\hat w_k^i \cdot z_i + \hat b_k^i) + c.\]

\begin{proposition} Given a compact set $K\subset\R^n$, for each $f\in F_1^{d_{\max}}$ there exist a Lyapunov function $V_f$ such that the following holds. For each $\eps>0$ the network depicted in and described after Figure \ref{fig:Lf_nn_l1} with $\sigma^1:\R\to\R$ infinitely differentiable and not polynomial, provides an approximation $\inf_{\theta\in\R^P} \|W(x;\theta) - V_f(x)\|_{\infty,K} \le \eps$ for all $f\in F_1^{d_{\max}}$ with a number of $N$ of neurons satisfying
\[ N = \OO\left(n^{d_{\max}+1}\eps^{-d_{\max}}\right).\]
\label{prop:V1}\end{proposition}

\vspace*{-4mm}

{\bf Proof:} Since the functions $V_i$ in Definition \ref{def:sg} are $C^1$ and the functions $\eta_i$ and $\alpha_i$ are continuous, the functions $\widetilde V_i$ from Theorem \ref{thm:sg} are $C^1$. We choose $\mu>0$ maximal such that $\mu \widehat V_i$ lies in $W_1^{d_i}$ and set $V_f=\mu V$ with $V=\sum_{i=1}^s \widehat V_i$ from Theorem \ref{thm:sg}.

Then, by Theorem \ref{thm:univ} we can conclude that there exist values $\hat a_k^i$, $\hat b_k^i$, $\hat w_k^i$, $c_i$, $k=1,\ldots,d_i$, such that the output of each sublayer $L_i$ satisfies
\[ \left\| \sum_{k=1}^{d_i} \hat a_k^i\sigma^1(\hat w_k^i \cdot z_i + \hat b_k^i) + c_i - \mu \widehat V_i\right\|_{\infty,K} \le \eps/n \]
for a number of neurons $M=\OO\left((\eps/n)^{-d_{\max}}\right)=\OO\left(n^{d_{\max}}\eps^{-d_{\max}}\right)$. 
Since this is true for all $L_1$, $\ldots$, $L_s$, by summing over the sublayers we obtain $W(x;\theta) = \sum_{i=1}^s \sum_{k=1}^{d_i} \hat a_k^i\sigma^1(\hat w_k^i \cdot z_i + \hat b_k^i) + c_i$ and thus
\[ \| W(x;\theta) - V_f(x)\|_{\infty,K} \le \eps \]
with the overall number of neurons $N \le nM = \OO\left(n^{d_{\max}+1}\eps^{-d_{\max}}\right)$.
\qed

We note that $\sigma^1$ is completely arbitrary under the stated assumptions. A systematic study of the performance of different $\sigma^1$ will be subject of future research.

\subsection{The case of unknown subsystems}

The approach in the previous subsection requires the knowledge of the subsystems $\Sigma_i$ in order to design the appropriate neural network. This is a rather unrealistic assumption that requires a lot of preliminary analysis in order to set up an appropriate network. Fortunately, there is a remedy for this, which moreover applies to a larger family of systems than $F_1^{d_{\max}}$ considered above. To this end, we consider the family $F_2^{d_{\max}}$ of maps $f$ in \eqref{eq:sys} with the following property: after a linear coordinate transformation $\tilde x=Tx$, $T\in\R^{n\times n}$ invertible and depending on $f$, the map $\tilde f(\tilde x) := Tf(T^{-1}\tilde x)$ lies in $F_1^{d_{\max}}$. In contrast to Section \ref{sec:known}, now we do {\em not} assume that we know the subsystems $\widetilde \Sigma_i$ of $\tilde f$, nor their dimensions $\tilde d_i$ and not even their number $\tilde s$.

\begin{figure}[htb]
\begin{center}
\includegraphics[width=8cm]{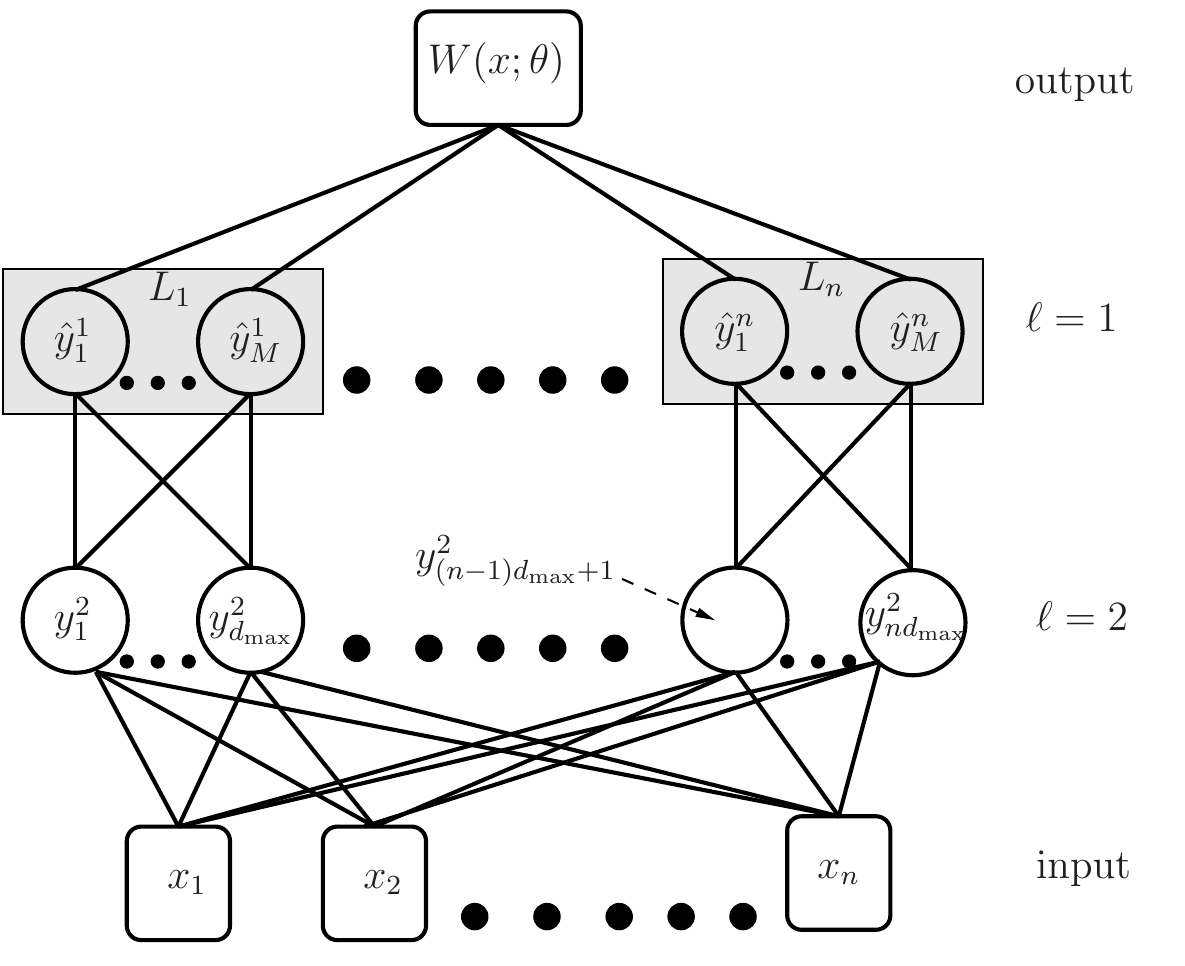}
\end{center}
\caption{Neural network for Lyapunov functions, $f\in F_2^{d_{\max}}$}
\label{fig:Lf_nn_l2}
\end{figure}

The neural network that we propose for $f\in F_2^{d_{\max}}$ is depicted in Figure \ref{fig:Lf_nn_l2}. 

Here, we use different activation functions $\sigma^\ell$ in the different levels. While in layer $\ell=1$ we use the same $\sigma^1$ as in Proposition \ref{prop:V1}, in Level $\ell=2$ we use the identity as activation function, i.e., $\sigma^2(x)=x$. Layer $\ell=1$ consists of $n$ sublayers $L_1$, $\ldots$, $L_n$, each of which has exactly $d_{\max}$ inputs and $M$ neurons. The coefficients and neuron values of each $L_i$ are again denoted with $\hat w_k^i$, $\hat a_k^i$, $\hat b_k^i$ and $\hat y_k^i$, respectively, for $k=1,\ldots,d_{\max}$. The $d_{\max}$-dimensional input of each neuron in $L_i$ is given by 
\[ (y^2_{(i-1)d_{\max}+1},\ldots,y^2_{i d_{\max}})^T =:\bar y^2_i.\] 
We note that this network is a special case of the lower network in Figure \ref{fig:gen_nn}. 

\begin{theorem} Given a compact set $K\subset\R^n$, for each $f\in F_2^{d_{\max}}$ there exist a Lyapunov function $V_f$ such that the following holds. For each $\eps>0$ the network depicted in and described after Figure \ref{fig:Lf_nn_l2} with $\sigma^1:\R\to\R$ infinitely differentiable and not polynomial in layer $\ell=1$ and $\sigma^2(x)=x$ in layer $\ell=2$, provides an approximation $\inf_{\theta\in\R^P} \|W(x;\theta) - V_f(x)\|_{\infty,K} \le \eps$ for all $f\in F_2^{d_{\max}}$ with a number of $N$ of neurons satisfying

\[ N = \OO\left(n d_{\max}+n^{d_{\max}+1}\eps^{-d_{\max}}\right).\]
\label{thm:V2}\end{theorem}
{\bf Proof:} Let $\tilde d_i$ be the (unknown) dimensions of the subsystems $\widetilde\Sigma_i$ and $p_i = 1+ \sum_{k=1}^{i-1} \tilde d_k$ the first index of the variables $\tilde z_i$ of $\widetilde \Sigma_i$, i.e., $\tilde z_i=(\tilde x_{p_i},\ldots, \tilde x_{p_{i+1}-1})^T$. Using the notation from above and the fact that $\sigma^2(x)=x$, the value of the inputs $y^2_k$ of the sublevels is given by 
\[ y^2_k = w_k^2 \cdot x + b_k^2. \]
Hence, by choosing $b_k^2=0$ and $w_k^2$ as the transpose of the \mbox{$j$-th} row of $T^{-1}$, we obtain $y^2_k = \tilde x_j$. Hence, by appropriately assigning all the $w_k^2$, we obtain 
\[ \bar y^2_i = \left(\begin{array}{c} \tilde z_i\\ 0 \\ \vdots \\ 0\end{array}\right),\]
where the number of the zeros equals $d_{\max}-d_i$. This can be done for $i=1,\ldots,\tilde s$, with $\tilde s$ being the number of subsystems. The inputs for the remaining sublayers $L_{\tilde s+1}$, $\ldots$, $L_n$ are set to $0$ by setting the corresponding $w^2_k$ and $b^2_k$ to $0$. For this choice of the parameters of the lower layer, each sublayer $L_i$ of the layer $\ell=1$ receives the transformed subsystem states $\tilde z_i$ (and a number of zeros) as input, or the input is $0$.

Since the additional zero-inputs do not affect the properties of the network, the upper part of the network, consisting of the hidden layer $\ell=1$ and the output, has exactly the structure of the network used in Proposition~\ref{prop:V1}. We can thus apply this proposition to the upper part of the network and obtain that it can realize a function $W(\tilde z;\theta)$ that approximates a Lyapunov function $\widetilde V$ for $\tilde f$ in the sense of Proposition \ref{prop:V1}. 

As the lower layer realizes the coordinate transformation $\tilde x = Tx$, the overall network $W(x;\theta)$ then approximates the function $V(x) := \widetilde V(Tx)$. By means of the invertibility of $T$ and the chain rule one easily checks that this is a Lyapunov function for $f$. The claim then follows since the number of neurons $N^1$ in the upper layer is equal to that given in Proposition \ref{prop:V1}, while that in the lower layer equals $N^2=nd_{\max}$. This leads to the overall number of neurons given in the theorem.
\qed

While we do not expect that more than two layers will improve Theorem \ref{thm:V2}, they may still perform better in practice. This will be subject of future research.

\section{Discussion}\label{sec:discussion}

In this section we discuss a few aspects and possible extensions of the results just presented.

\begin{enumerate}[(i)]
\item From the expressions for $N$ in Proposition \ref{prop:V1} and Theorem \ref{thm:V2} one sees that for a given $\eps>0$ the storage effort only grows polynomially in the state dimension $n$, where the exponent is determined by the maximal dimension of the subsystems $d_{\max}$. The proposed approach hence avoids the curse of dimensionality, i.e., the exponential growth of the effort. There is, however, an exponential dependence on the maximal dimension $d_{\max}$ of the subsystems $\Sigma_i$ in the small-gain formulation. This is to be expected, because the construction relies on the low-dimensionality of the $\Sigma_i$ and if this is no longer given, we cannot expect the method to work efficiently.

\item In order to use the proposed architecture in practice, suitable learning methods have to be developed. This is a serious problem in its own right that is not addressed in this paper. However, given that the decisive property of a Lyapunov function \eqref{eq:lf} is a partial differential inequality, there is hope that similar methods as for the learning of solutions of partial differential equations, see, e.g., \cite{HaJE18,HuPW19,SirS18}, can be implemented. This will be subject of future research.

\item There have been attempts to use small-gain theorems for grid-based constructions of Lyapunov functions, e.g., in \cite{CaGW09,Li15}. The problem of this construction, however, is, that it represents the functions $\hat V_i$ from Theorem \ref{thm:sg} separately for the subsystems and the small-gain condition has to be checked additionally (which is a difficult task). The representation via the neural network does not require to check the small-gain condition nor is the precise knowledge of the subsystems necessary. 

\item \label{it:u} The reasoning in the proofs remains valid if we replace $f(x)$ by $f(x,u)$ and asymptotic stability with ISS. Hence, the proposed network is also capable of storing ISS and iISS Lyapunov functions.

\item \label{it:relu} In current neural network applications ReLU activation functions $\sigma(r) = \max\{r,0\}$ are often preferred over $C^\infty$ activation functions. The obvious disadvantage of this concept is that the resulting function $W(x;\theta)$ is nonsmooth in $x$, which implies the need to use concepts of nonsmooth analysis for interpreting it as a Lyapunov function. While one may circumvent the need to compute the derivative of $W$ by means of using nonsmooth analysis or by passing to an integral representation of \eqref{eq:lf} (see also item \eqref{it:C0}, below), the nonsmoothness may cause problems in learning schemes and it remains to be explored if these are compensated by the advantages of ReLU activation functions.

\item \label{it:maxlf} There are other types of Lyapunov function constructions based on small-gain conditions different from Definition \ref{def:sg}, e.g., a construction of the form 
\[ V(x) = \max_{i=1,\ldots,s} \rho_i^{-1}(V_i(z_i)), \]
found in \cite{Ruef07,DaRW10}. Since maximization can also be efficiently implemented in neural networks (via so-called max pooling in convolutional networks), we expect that such Lyapunov functions also admit an efficient approximation via deep neural networks. However, when using this formulation we have to cope with two sources of nondifferentiability that complicate the analysis. On the one hand, this is caused by the maximization in the definition of $V$ and on the other hand by the functions $\rho_i^{-1}\in\KK_\infty$, which in most references are only ensured to be Lipschitz.

\item \label{it:C0} One may argue that a mere approximation of $V(x)$ is not sufficient for obtaining a proper approximation of a Lyapunov function, because in order for $W(x;\theta)$ to be a Lyapunov function, $d/dx W(x; \theta)$ should also approximate $DV(x)$ such that \eqref{eq:lf} is maintained approximately.

On the one hand, this is not strictly true. It is known (see, e.g., Section 2.4 in \cite{GruP17}), that under mild regularity conditions on $f$  asymptotic stability follows from the existence of a $T>0$ with
\begin{equation} V(x(T,x_0)) \le V(x_0) - \int_0^T h(x(t,x_0)) dt \label{eq:intLf}\end{equation}
for all $x_0\in \R^n$. This integral inequality, in turn, holds for each Lyapunov function simply by integrating \eqref{eq:lf}. Now, in order to ensure that a function $W$ satisfies \eqref{eq:intLf} approximately (i.e., plus an error term $\eps$ on the right hand side) it is sufficient to approximate $V$ itself by $W$. An approximation $DW\approx DV$ is not needed. 

On the other hand, it is known that neural networks also approximate derivatives of functions, cf.\ \cite{HoSW90}. Hence, it seems likely that our results are extendable to an approximation of $V$ and $DV$. This will be investigated in future research.
\end{enumerate}

\section{Conclusion}
We have proposed a class of deep neural networks that allows for approximating Lyapunov functions for systems satisfying a small-gain condition. The number of neurons needed for an approximation with fixed accuracy depends exponentially on the maximal dimension of the subsystems in the small-gain condition but only polynomially on the overall state dimension. Thus, it avoids the curse of dimensionality. The network structure does not need any knowledge about the exact dimensions of the subsystems and even allows for a subsystem structure that only becomes visible after a linear coordinate transformation.

The result suggests that suitably designed deep neural networks are a promising efficient approximation architecture for computing Lyapunov functions. It motivates future research on the development of suitable training techniques in order to actually compute the Lyapunov functions.

%\begin{ack}
%Place acknowledgments here.
%\end{ack}

\bibliography{../../bib/lars_consolidated}             % bib file to produce the bibliography
                                                     % with bibtex (preferred)
                                                   
%\begin{thebibliography}{xx}  % you can also add the bibliography by hand

%\bibitem[Able(1956)]{Abl:56}
%B.C. Able.
%\newblock Nucleic acid content of microscope.
%\newblock \emph{Nature}, 135:\penalty0 7--9, 1956.

%\bibitem[Able et~al.(1954)Able, Tagg, and Rush]{AbTaRu:54}
%B.C. Able, R.A. Tagg, and M.~Rush.
%\newblock Enzyme-catalyzed cellular transanimations.
%\newblock In A.F. Round, editor, \emph{Advances in Enzymology}, volume~2, pages
%  125--247. Academic Press, New York, 3rd edition, 1954.

%\bibitem[Keohane(1958)]{Keo:58}
%\newblock \emph{Power and Interdependence: World Politics in Transitions}.
%\newblock Little, Brown \& Co., Boston, 1958.

%\bibitem[Powers(1985)]{Pow:85}
%\newblock Is there a way out?
%\newblock \emph{Harpers}, pages 35--47, June 1985.

%\bibitem[Soukhanov(1992)]{Heritage:92}
%\newblock \emph{{The American Heritage. Dictionary of the American Language}}.
%\newblock Houghton Mifflin Company, 1992.

%\end{thebibliography}

                                                                         % in the appendices.
\end{document}